\documentclass[a4paper,10pt]{article}
\usepackage{amsfonts}
\usepackage{amssymb}
\usepackage[utf8]{inputenc}

\usepackage[english ]{babel}

\usepackage{bbding}

\usepackage[all]{xy}

\usepackage{tikz-cd} 
\usepackage[all]{xy}
\usepackage{amsfonts}
\usepackage{amssymb}
\usepackage{graphicx}
\usepackage{mathtools}
\usepackage{skak} 
\usepackage{foekfont}

\usepackage{calligra}
\usepackage[T1]{fontenc}

\usepackage{amsmath,mathdots}

\usepackage{stmaryrd}

\usepackage[utf8]{inputenc}

\usepackage{mathtools}
\usepackage{mathdots}
\usepackage{tikz}

\definecolor{amber(sae/ece)}{rgb}{1.0, 0.49, 0.0}
\newfont{\rsfsten}{rsfs10 scaled 1200}
\usepackage{MnSymbol}
\usepackage{tikz}
\usepackage{amscd}
\usepackage{MnSymbol}
\usepackage{wasysym}

\usepackage{mathrsfs}
\usepackage{pdfcolmk}

\usepackage{graphicx}
\newcommand*{\rom}[1]{\expandafter\@slowromancap\romannumeral #1@}

\usepackage{wrapfig}

\newcommand{\tightunderset}[2]{%
  \mathop{#2}\limits_{\vbox to .3ex{\kern-0.95ex\hbox{$#1$}\vss}}}
\newcommand{\tightoverset}[2]
{%
  \mathop{#2}\limits_{\vbox to .3ex{\kern-0.95ex\hbox{$#1$}\vss}}}

\newcommand{\oset}[2]{%
  {\mathop{#2}\limits^{\vbox to -.5\ex@{\kern-\tw@\ex@
   \hbox{\scriptsize #1}\vss}}}}
\usepackage[utf8]{inputenc}

\usepackage {ifsym}

\usepackage {url}

\title {No metrics with  Positive   Scalar Curvatures on Aspherical 5-Manifolds }

\author{Misha Gromov}

\begin{document}
\maketitle
\tableofcontents
\begin{abstract} 

A metric space $X$ is called {\it uniformly acyclic} if there there exists an {\it acyclicty control    function} $R=R(r)=R_X(r)\geq r $, $0\leq r  <\infty$,  such that the homology inclusion homomorphisms between the balls around all points $x\in X$,
$$H_i(B_x(r))\to H_i(B_x(R))$$
{\it vanish} for all $i=1,2,... $ .

 We show that if  
 
 {\sf a complete orientable  $m$-dimensional manifold $\tilde X$ of dimension $m\leq 5$  admits a proper (infinity goes to infinity)  distance decreasing  map to a complete  $m$-dimensional { \it uniformly acyclic} manifold, then the scalar curvature of $\tilde X$ {\it can't be uniformly positive}, 
 $$\inf _{x\in \tilde X}Sc(X,x) \leq 0.$$}

Since the universal coverings $\tilde X$  of compact aspherical manifolds  $X$  are  {\it uniformly acyclic}, (in fact,  {\it uniformly contractible}), these $X$, admit no metrics with $Sc>0$  for $dim (X)\leq 5$.

Our argument, that depends on {\it torical symmetrization} of  {\it stable $\mu$-bubbles}, is  inspired by  the recent paper [C\&L 2020]  by Otis Chodosh  and  Chao Li on non-existence of metrics with $Sc>0$  on aspherical 4-manifolds  and is also  influenced by the ideas 
 of Jintian Zhu and Thomas Richard 
from  their papers [J.Z. 2019] and [T.R. 2020].

 %On the 2-systole of stretched enough positive ...
%arxiv.org ? math
%Jul 6, 2020 - On the 2-systole of stretched enough positive scalar curvature metrics on S 2 x S 2. Authors:Thomas Richard (LAMA). Download PDF. Abstract: ...
%by T Richard - ?2020
\end{abstract}

\section{$ \mathbb T^\ast$-Stabilization of   Scalar Curvature} 

Let $X=(X,g(x))$ be a Riemannian manifold, let  $\Phi$ be an $N$-tuple of  positive smooth functions   $\phi_i(x)$, $i=1,2,..., N$, on $X$, and let 
$$g^{\ast}=g^{\ast}_\Phi=g(x)+\phi^2_1(x)dt_1^2+\phi^2_2(x)dt_2^2+...+\phi^2_N(x)dt_N^2$$
be the iterated warped product metric on $X^\ast =X\times \mathbb R^N$.

Observe that 

$\bullet_1$ {\sf the metric $g^{\ast}$ is $\mathbb R^N$-invariant under the natural action of  $\mathbb R^N$ on $X^\ast$, where $X$ is identified with the quotient space $X^\ast/\mathbb R^N$  for this action;

$\bullet_2$ the natural  (quotient) map $X^\ast\to X$ is a {\it Riemannin submersion}:  it is {\it isometric} on the smooth {\it horizontal  
curves} $C\subset  X^\ast$, i.e. those normal to the  {\it $\mathbb R^N$-fibers} that are 
$\mathbb R^N_x\subset X^\ast, x\in X$;

$\bullet_3$ the zero-section embedding $X\to X^\ast$ is an {\it isometry. }\hspace {1mm}}

Furthermore, a straightforward computation (compare   [F-C\&S 1980], [G\&L  1983], [J.Z. 2019]) %D. Fisher-Colbrie & R. Schoen, The structure of complete stable minimal
%surfaces in 3-manifolds of non-negative scalar curvature, Comm. Pure
%Appl. Math., 33 (1980) 199-211.
%Mikhael Gromov and H. Blaine Lawson, Jr., Spin and scalar curvature in the presence
%of a fundamental group. I, Ann. of Math. (2) 111 (1980), no. 2, 209{230. MR 569070
%[6] , Positive scalar curvature and the Dirac operator on complete Riemann-
%ian manifolds, Inst. Hautes Etudes Sci. Publ. Math. (1983), no. 58, 83{196 (1984).\S12
%section 2.4 in Four lectures on scalar curvature, arXiv:1908.10612
%Rigidity of Area-Minimizing 2-Spheres in n-Manifolds with Positive Scalar Curvature. by J Zhu  	arXiv:1903.05785
shows that\hspace {1mm}

$\bullet_4$ {\sf   the scalar curvature of $X^\ast$, which,   being   $\mathbb R^N$-invariant,  is regarded as a function on $X$, satisfies:
$$Sc(g^\ast,x)=Sc(g, x) -2\sum_{i=1}^N\frac{\Delta_g(\phi_i(x))}{\phi_i(x)}  -2\sum _{i<j} \langle\nabla_g(\log \phi_i(x), \nabla_g(\log \phi_i(x)\rangle.$$}

Manifolds  $X^\ast$, which are  defined  as above     by $N$-tuples $\Phi$  are called 
{\it warped  $\mathbb T^N$-extensions of $X$.}

{\it From  $X^\ast$ to $\underline X^\ast =X^\ast/\mathbb Z^N$.} In what follows, we prefer to work with  $X^\ast$ divided  by the action of the lattice $\mathbb Z^N\subset \mathbb R^n$, where  the Euclidean fibers $\mathbb R^N_x\subset X^\ast$ become  torical, that are
$\mathbb T_x^n=  \mathbb R_x^n/\mathbb Z^N$; this explains our  
$\mathbb T^\ast$-terminology.

These $X^\ast$ and $\underline X^\ast$ are involved in our  arguments  in two ways,
\vspace{1mm}

I. {\sc $\mathbb T^\ast$-Stabilization Principle}:  {\it Many geometric properties of manifolds $X$ implied by the inequality  $Sc(X)\geq \sigma$ follow (possibly in a weaker form)  from the inequality $Sc(X^\ast)\geq \sigma(x)$  satisfied by some     warped $\mathbb T^N$-extension $X^\ast$  of $X$.}

 \vspace {1mm}

 {\it Example 1: Weakened  $\mathbb T^\ast$-Stable 2d Bonnet-Myers Diameter Inequality}. {\sf If a closed  connected surface $X$ admits 
a  warped $\mathbb T^N$-extension $X^\ast$ with 
   $Sc^\ast(X)\geq \sigma$, then  the  diameter of $X$ is  bounded as follows.
 $$diam (X)\leq  2\pi\sqrt{\frac { N+1}{ (N+2)\sigma}} \hspace {1mm}< 
 \frac {2\pi}{\sqrt\sigma}.  $$} \vspace {1mm}

{\it Proof.}  Given two points $x_1,x_2 \in X$, take two small $\varepsilon$-circles $Y_{-1}$  and $Y_{+1}$  around them, let $X_\varepsilon \subset X$ be the band between them and 
 and apply the (elementary) $\mathbb T^N$-invariant  case of  the (over)-{\it torical  $\frac {2\pi}{n}$-Inequality} from [M.G. 2018].
%	arXiv:1710.04655 Metric Inequalities with Scalar Curvature. Authors:Misha Gromov

{\it Remarks.} (a)   This  proof is similar to that of  theorem 10.2 in  [G\&L 1983],  where the latter
concerns stable minimal surfaces $Y$    in 3-manifolds $X$ with $Sc(X)\geq \sigma>0$,  
 and  delivers the following bound on the filling radii of contractible (actually homologous to zero) curves $S\subset X$:
$$fil.rad(S\subset X) \leq  \frac {2\pi}{\sqrt\sigma},$$ 
and where
 a version of our  example  1 for $N=1$ is implicit   in the calculations on pp 178-180.

(The above proof gives a slightly better bound in this case,  namely $fil.rad(S\subset X) \leq  \frac {\pi\sqrt 2}{\sqrt\sigma}$.)

(b) It is not hard to adapt   the warped product metric on $ \mathbb T^{n-1}\times [-1, +1]$ that was  used for showing  optimality of the  $\frac{2\pi}{n}$ -inequality  (see section 2 in [M.G. 2018]),  for proving the optimality of the  above inequality  for closed surfaces.

 (c) If $X$ is a complete  Riemannian  3-manifold with $H_1(X;\mathbb Q)=0$, e.g. homeomorphic to the sphere $S^3$ or to a connected sum of space forms, the above argument shows that the inequality $Sc(X^\ast)\geq \sigma>0$ for a warped $\mathbb T^N$-extension  of $X$ implies that the distances $d$ between pairs of circles in $X$ with {\it non-zero} linking numbers are bounded by $2\pi\sqrt{(N+2)/(N+3)\sigma}$.
Thus, for instance, if $(S^3,g)\to (S^2,h)$ is a {\it non-contractible distance 
non-increasing} map and if  $Sc(g)\geq \sigma>0$, then 
$diam_h(S^2)\leq 2\pi\sqrt{2/3\sigma}$.

 \vspace {1mm}
 
  {\it Example 2: Sharp Stabilization of the Gauss-Bonnet inequality.}   {\sf If a closed  surface $X$ admits 
 a warped $\mathbb T^N$-extension $X^\ast$ with 
   $Sc^\ast(X)\geq 2$, then  the area of $X$ satisfies
   $$ area(X)\leq 4\pi.$$}
   
  \vspace {1mm}

{\it Commentary.} This inequality,  proven by Jintian Zhu in [J.Z. 2019], %Rigidity of Area-Minimizing 2-Spheres in n-Manifolds with Positive Scalar Curvature. by J Zhu  	arXiv:1903.05785, 
is not used in the sequel; we brought it up, since it  served as a  motivation for  the diameter  inequality in  our Example 1. 

II {\it Torical Symmetrization.}  {\sf Let $X$ be a closed orientable Riemannin manifold of dimension $n+N$  an  let 
$$f: X\to \mathbb T^N$$
be a continuous map. }  

{\it Theorem} 1.  If $dim(X)\leq 8$,  then\vspace{1mm}

\hspace {-6mm}{\it there exists

$\bullet$ a smooth  closed (possibly disconnected)  orientable submanifold
 $Y\subset X$ of dimension $n$ in the homology class of the pullback $f^{-1}(t)$, $t\in \mathbb T^N$,  that is, in other words, the Poincar\'e dual of the pullback $f^\ast [\mathbb T^n]\in H^N(X)$  of the fundamental cohomology class  $[\mathbb T^N]\in H^N( \mathbb T^N)$,}  
 and 

$\bullet $ {\it a warped $\mathbb T^N$-extension $Y^\ast $ of $Y$ such that  the scalar curvature of $Y^\ast$ is bounded from below by the scalar curvature of $X$ as follows.
$$\mbox {$Sc(Y^\ast,y)\geq Sc(X,y)$   for all $y\in Y\subset X$.}$$}
\vspace{1mm}

{\it Commentaries.}  (a) This kind of  symmetrization for $n=2$ and $N=1$  goes back to the 1980-paper by Fisher-Colbrie and  Schoen    [F-C\&S 1980],  %The structure of complete stable minimal
%surfaces in 3-manifolds of non-negative scalar curvature, Comm. Pure
%Appl. Math., 33 (1980) 199-211.
 which was extended   to other  dimensions $n$ and $N$  in [G\&L 1983]  with a  use of the  Schoen-Yau descent method by minimal hypersurfaces %Mikhael Gromov and H. Blaine Lawson, Jr., Spin and scalar curvature in the presence
%of a fundamental group. I, Ann. of Math. (2) 111 (1980), no. 2, 209{230. MR 569070
%[6] , Positive scalar curvature and the Dirac operator on complete Riemann-
%ian manifolds, Inst. Hautes Etudes Sci. Publ. Math. (1983), no. 58, 83{196 (1984).\S12
  and    applied  in [M.G. 2018] to manifolds with boundaries.%arXiv:1710.04655 Metric Inequalities with Scalar Curvature. Authors:Misha Gromov

(b)  Our  "submanifold"  $Y$ is,  in fact, a disjoint union of several components that may come with {\it integer multiplicities}  $\neq \pm1$,  but  this is irrelevant for our applications. Besides, we limit 
ourselves to  those situations, where   one can  slightly move 
 identical copies  of these components {\it without affecting our inequalities}, by  replacing 
$kY$ by $\cup_{i=1}^kY_i$ with  $Y_i$ being small normal  mutually disjoint shifts of $Y$. 

\section {Stable $\mu$-Bubbles in Riemannian Bands} 

 Let $X$ be a compact   {\it Riemannian band} (in the sense of [M.G. 2018]),  of dimension $m$   (also called {\it condenser} in [M.G. 2019]), %Four lectures on scalar curvature, arXiv:1908.10612 
 that is a compact Riemannin manifold with a boundary, where this boundary $\partial X$   is decomposed into two non-empty mutually disjoint  parts,
$$\partial(X)=\partial_-\cup \partial_+,$$
and where both parts $\partial_-\subset \partial X$ and $ \partial_+\subset \partial X$  are unions of  connected components of  $\partial X$.

{\sf Let $d_\pm=d_\pm(X)$
denote the {\it width of $X$}, that is  
 $$d_\pm=width(X)=dist(\partial_-,\partial_+)$$
 and let the scalar curvature of $X$ be bounded from below by $\sigma>0$,
$$Sc(X)\geq \sigma.$$}

{\it Theorem} 2. {\sf Let 
$$\gamma=\sigma d_\pm^2\geq \frac {4\pi^2(m-1)}{m}.  $$}
If $m=dim(X)\leq 8$, then  

\hspace {-6mm}{\it there exists a smooth hypersurface $Y\subset X$,  which separates $\partial_-$ from $\partial_+$, 
 and a  warped 
 $\mathbb T^1$-extension $Y^\ast$ of $Y$,   such that
 $$Sc(Y^\ast)\geq \beta_m\sigma,$$  
where $\beta_m=\beta_m(\gamma)$ is a       monotone increasing  function in  $\gamma $ for 
$  \frac {4\pi^2(m-1)}{m}\leq \gamma<\infty$, such that
$$ \beta_m\left (\frac {4\pi^2(m-1)}{m}\right)=0,$$
$$\beta_m(\gamma)\to 1\mbox { for } \gamma \to \infty,$$
 $$\beta_m(\gamma) >0 \mbox { for } \gamma> \frac {4\pi^2(m-1)}{m};$$ 
 moreover,
 $$\beta_m(m(m-1)\pi^2)\geq \frac {m-2}{m}.\leqno {[\Circle^m]}$$}

{\it About the Proof.}  The required separating hypersurface $Y\subset X$ is obtained as a  minimizing hypersurface for  
 the functional
 $$Y\mapsto vol_{m-1}(Y)- \int_{X_-}\mu(x)dx,$$
for a suitable 
function $\mu(x)$ on $X$,  where $X_-\subset X$ denotes the  band between $Y$ and $\partial_-$,
see section 5.4 in [M.G. 2019], where the inequality $[\Circle^m]$ follows by a comparison argument applied to
 the unit  sphere $S^m$ minus two opposite points and where   the spherical metric $g_m$ on $S^m$ is  decomposed in polar coordinates as a   warped product  
 of  $S^{m-1}$ with   $(-\pi/2, \pi/2)$,
 $$ g_m=(\cos t )^2 g_{m-1} +dt^2.$$

\vspace {1mm}

  Torical symmetrization (II from  the previous  section) applied to the above $\mu$-bubble  yields the following.

{\it Corollary.} {\sf Let a  compact  oriented Riemannin band $X$  of dimension $m= n+N+1$ admit a continuous  {\it band map} $f$  to $\mathbb T^{N}\times [-1,1]$, (i.e.
 $f:\partial_\pm\to \mathbb T^{N}\times \{\pm1\}$).

  Let $$Sc(X)\geq \sigma>0$$
and 
 $$\sigma d_\pm^2\geq m(m-1)\pi^2,$$
where  $d_\pm$  denotes the width of the band $X$ as earlier.}

If $m\leq 8$, then 

{\it the homology class of the $f$-pullback of a point  $t\in \mathbb T^{N}\times [-1,1]$} (i.e. dual to 
$ f^\ast  [\mathbb T^{N}\times [-1,1]]\in H^{N+1}(X, \partial X)$ of the 
 the  relative fundamental   class of $\mathbb T^{N}\times [-1,1]$)
{\it can be realized by a  closed smooth $n$-dimensional submanifold $Y\subset X$, which admits a warped  
$\mathbb T^{N+1}$-extension  $Y^\ast$,
such that  
$$Sc(Y^\ast)\geq \frac  {m-2}{m}\sigma.$$}

{\it Commentaries} (a) The above torical   symmetrization  of $\mu$-bubbles is  similar to that
used by Thomas Richard in  [T.R. 2020] %On the 2-systole of stretched enough positive scalar curvature metrics 
%on $S^2 \times  S^2$. arXiv:2007.02705
used for  {\sl a bound on the 2-systole  of metrics with $Sc \geq \sigma>0$  on $S^2\times  S^2$}.

(b) Probably, the desingularization result  for minimizing hypersurfaces  by Joachim Lohkamp [J.L. 2018]  %J. Lohkamp [Loh(smoothing) 2018] {\sl Minimal Smoothings of Area Minimizing Cones}, https://arxiv.org/abs/1810.03157
  would allow the above to hold for $m\geq 9$ as well.

(c) Probably    comparison of  $X$ with suitable  non-spherical warped products (see 5.4 in [M.G. 2019]) would lead 
to  a better estimate of  the function  $\beta_m$ and  allow sharpening of the 
  above inequalities.

\section {$\mu$-Bubbles in Codimension 2}

The following proposition, that is adapted from [R.T. 2020] (compare [J.Z. 2019]),  plays a key role in the proof of our main result.

{\it Richard's Lemma.} {\sf Let $X$ be an   oriented $m$-dimensional Riemannin   manifold (possibly non-compact and non-complete)  with compact  boundary 
and $X_0\subset X$  be an open subset  with smooth boundary  such that the complement $X\setminus X_0$ is compact.
Let $h\in H_{m-2}( \partial X)$ and $h_0\in H_{m-2}(X_0)$ be homology classes, which have {\it equal  images}  under  the  homomorphisms      induced by the inclusions\ $\partial X\hookrightarrow X \hookleftarrow X_0$,
that are
$$h\in H_{m-2}(\partial X)\rightarrow   H_{m-2}(\partial X) \leftarrow H_{m-2}(X_0)\ni h_0.$$}

{\it Let   $$Sc(X)\geq \sigma>0, $$ 
and
$$ dist^2 (X_0, \partial X)  \geq \frac  {m(m-1)\pi^2}{\sigma}.$$}

If $m\leq 8$, then 

{\it the image of the homology class $h$ in $H_{m-2}(X)$ can be realized by a  closed smooth $(m-2)$-dimensional submanifold $Y\subset X$, which admits a warped  
$\mathbb T^{2}$-extension  $Y^\ast$,
such that  
$$Sc(Y^\ast)\geq \frac  {m-2}{m}\sigma.$$}  
 
 {\it Proof.} Let $ \overleftrightarrow X =X\setminus X_0\subset X$ be the band between $\partial_-=\partial X_0$ and   $\partial_+=\partial X$  and  let  $ \overleftrightarrow h\in H_{m-1} ( \overleftrightarrow X, \partial \overleftrightarrow X)$ be the 
 relative   class that establish homology equivalence  between $h\in H_{m-2}(\partial X)$ and 
 $h_0\in H_{m-2}(X_0)$, where the latter is  moved from   $H_{m-2}(X_0)$  to $H_{m-2}(\partial X_0)$  by the excision property.

  Let $\overleftrightarrow  h^\perp  \in H^1(\overleftrightarrow X)$ be the integer cohomology class Poincar\'e dual to   $\overleftrightarrow h$ and let  us 
  induce $\overleftrightarrow h^\perp$ from the fundamental class of the circle by a continuous map 
  $f_\perp:\overleftrightarrow X\to \mathbb T^1$.
 
 Then, using the distance function  $x\mapsto dist (x,X_0)$, we  construct a band map 
 $\ \overleftrightarrow f:\overleftrightarrow X\to [-1,1] $ and apply  the above corollary
 to the band map defined by the pair $(f_\perp, \overleftrightarrow f)$,
 $$(f_\perp,  \overleftrightarrow f): \overleftrightarrow X \to  \mathbb T^1\times  [-1,1]. $$
  
 %${\pigpenfont N}$
  %$\Ds$
 
%$\diamonddot$

 \section {Recollections on Filling Radius, Uryson Width and Uniform Acyclicty}

%Katz, M. Systolic geometry and topology. With an appendix by J. Solomon. Mathematical Surveys and Monographs, volume 137. American Mathematical Society, 2007.

%Gromov, M.: Filling Riemannian manifolds, J. Diff. Geom. 18 (1983), 1?147.

%Guth, Larry. ?Volumes of Balls in Riemannian Manifolds and Uryson Width.? Journal of Topology and Analysis (February 22, 2016)

%Guth, Larry. Metaphors in systolic geometry. Proceedings of the International
%Congress of Mathematicians 2010 (ICM 2010)

%Sep 9, 2019 - Uryson width and volume. Authors:Panos Papasoglu
%Geometric and Functional Analysis volume 30, pages574?587(2020)

%A short proof of Gromov's filling inequality

 The {\it  filling radius} $fil.rad(X)$,   also called the  {\it absolute filling radius} of an $n$-dimensional  orientable  manifold (or a pseudomanifold)  $X$ with a metric $dist_X$ on it is defined as the supremum of the numbers $r>0$, such that if  $X_+ \supset  X$ is a  {\it metric extension} of $X$, that is a  metric space that contains   $X$ and such that $dist_{X_+}|X=dist_X$, the fundamental class $[X]\in H_n(X)$ {\it doesn't vanish} (i.e. $X$ doesn't bound) in the $r$-neighbourhood $U_r(X)\subset X_+$.\footnote  {See [M.G. 1983], [M.K. 2007], [S.W 2007], [L.G. 2010],  [L.G. 2016], [P.P. 2020] for details and other properties   of the filling radius and other concepts discussed in this section.}

{\it Filling with Coefficients.} The above definition   makes sense for the fundamental class of $X$ in the homology with a given coefficient ring.
For instance one may speak of the  {\it rational filling radius} $fil.rad_\mathbb Q (X)$
which may be smaller than $fil.rad=fil.rad_\mathbb Z (X)$.

Moreover, if $\hat X\to X$ is finite covering map, then $fil.rad_\mathbb Q (X)\leq  fil.rad_\mathbb Z (\hat X)$.

(It seems unclear, in general.  what {\sf possible values of the ratio} $\frac {fil.rad(X)}{fil.rad(\hat X)}$ could be 
for various  Riemannin metrics on $X$ and  those on  $\hat X$ induced by the covering
$\hat X\to X$.)

4.A. {\it Cubic Example} If  $X$   can be  covered by $2n +2$ subsets  $Y_{\pm i}\subset X$, $i=0,1,...,n,$  
such that $dist_X(Y_{-i},Y_{+i})\geq d$, then
$fil.rad_\mathbb Q(X)\geq d/2$.

Indeed,  map $X_+$ to 
$\mathbb R^n_+$ by  $n$ distance  functions $x\mapsto dist_X(x,Y{_-i})$ and compose it with the radial projection from $\mathbb R^n_+$ to the $d$-cube $[0,d]^n\subset \mathbb R^n_+$.

The restriction of this composed map $F: X_+\to [0,d]^n$ sends $X$ to the boundary sphere 
$\partial [0,d]^n$ with degree 1.

Hence,  if an  $(n+1)$-chain $C \subset X_+$  bounds $X$, it is  is sent by $F$ {\it onto } this cube 
where  all   points $x_+\in C$ from the  $F$-pullback of the center of the cube, 
$(\tightunderset {n} {\underbrace {d/2,d/2,....,d/2}})\in  [0,d]^n$
have $dist(x_+, X)\geq d/2$. QED.

\vspace {1mm}

4.B*. {\it Corollary.}\footnote{This corollary and everything else marked with "*" is not used in the proof of our main theorems.} {\sf Let $X$ be a Riemannian manifold, let $Y\subset X$ be the shortest non-contractible closed curve and let $\hat Y\to Y$ be a finite covering    of this $Y$ (where the curve $\hat Y$ may be contractible in $X$).

Then the filling radius of $\hat Y$ for the metric on $\hat Y$ induced by distance function on $X$\footnote{Never mind vanishing of this metric on the pullbacks of he points by the covering map $\hat X\to X$.} is bounded by the length $l$ of $Y$ as follows
$$ fil.rad (\hat Y)\geq \frac {l}{8}.$$
Moreover, the same remains true for the shortest curve {\it non-homologous to zero}.
or more generally, being non-trivial in some quotient group of the fundamental group 
$\pi_1(X)$.}  

{\it Proof.} Divide $Y$ in  a square like fashion into  four  segments 
$s_{\pm i}$,         $i=1,2$, of lengths $\frac{l}{4}$.  

If $Y$ is the shortest, then the distances $d$ between both pairs of the  opposite segments must be 
$\geq \frac {l}{4}$; otherwise  $Y$ could be decomposed  into two curves of curves of lengths
$l_1+d$ and $l_2+d$, where $l_1,l_2\leq \frac {3l}{4}$ (and $l_1+l_2 =l$). QED.

(This, as  shall see in section 6,  implies  that closed orientable   3-manifolds $X$ with $Sc(X)\geq \sigma>0$  and   the  fundamental groups of which are {\it non-free}, contain closed non-contractible geodesics  of length 
$\leq 100/\sqrt\sigma$.) 
\vspace {1mm}

 {\it Uryson's $k$-width}, denoted  $width_k(X)$, of a metric space $X$ is the infimum of the  numbers $w$, such that $X$ admits a covering of multiplicity $\leq k+1$  by closed subsets $V_i\subset X$ with $diam(V_i)\leq w$.

Equivalently, -- assuming   that $X$ is compact  for   safety sake -- this is equal to  the infimum of $w$, such that $X$ admits a continuous map 
to a $k$-dimensional polyhedral space, 
$$f:X\to P,$$
such that the diameters of the $f$-pullbacks  of all points $p\in P$ satisfy
$$diam(f^{-1}(p))\leq w.$$

Indeed, given a covering  of $ X$ by $V_i$, let $P$ be the nerve of this covering  and let  $f:X\to P$ be the standard  continuous map defined via a partition of unity subordinated  to the    covering of $ X$ by small 
$\epsilon$-neighbourhoods $U_i\supset V_i$.

Conversely,  given a  continuous map $f\to P$, $dim(P)=k$,  take a covering of $P$ by sufficiently small subsets with multiplicity $k+1$ and then pull it back by $f$ to $X$.

\vspace{1mm}

4.C. {\sf The filling radius is bounded by the Uryson width as follows.}
$$fill.rad (X)\leq \frac{1}{2}width_k(X)\mbox { for  all  } k<dim (X).$$ 

Indeed, let  $ X_+$ be  the {\it cylinder $C_f$}  of the map $f$,   obtained by attaching the 
product $X\times [0,w/2]$ to $P$ by the map $f: X=X\times \{w/2\}\to P$
and observe that

$\bullet$  there is an extension of the metric of $X$ to $C_f$ that  keeps
the lengths of the  segments  $ [\{x\}\times [0,w/2]$, $x\in X$, equal to $ w/2.$

and that

$\bullet$ if $dim(P)<dim(X)$, then  $X$ is homologous to zero in $C_f\supset X= X\times \{0\}$.      

\vspace {1mm}

4.D. {\it Elementary Metric Lemma.} {\sf Let   $X$    be   a path metric space , e.g. a Riemannian manifold, possibly non-complete and or with a boundary and   let 
$\phi(x)$ be the distance function to a point $x_0\in X$,
$$\phi(x)=dist(x,x_0).$$}

{\sl If a connected component of a level of $\phi$  has diameter $\geq d$, 
%$L\subset A$  be a path connected subset 
%such that the distance function to a point $a_0\in A$ is constant on $L$, say 
%$$dist(a,a_0)=r\mbox  { for all $a\in L$}.$$} 
then  the space  $X$ contains a closed curve $S$, such that,  this curve $S$ itself and all its finite coverings $\hat S\to S$ have their filling radii bounded from below  by 
 $$fil.rad(\hat S)\geq   \frac {diam(L)}{6}.  $$}

{\it Proof.} Let $x_1, x_2\in \phi_{-1}(r)\in X$ 
be two point in a  connected component of an   $r$-level of $\phi$ for some $r>0$,
such that $dist(x_1,x_2)\geq d$, assume without loss of generality that this component is path connected and proceed as follows.

$\bullet$ Join   the points  $x_1$ and $x_2$    by a curve $l$ in this $r$-level of $\phi$,

%$\bullet$ let  $x_1$ and $x'_2\$
%be the nearest points to $x_1$ and to $x_2$  in $X_0$ join them by a curve $l'\subset X_0$,

$\bullet$  join $x_1$ and $x_2$  with $x_0$   by shortest segments  $s_1$  and $s_2$ in $X$,

$\bullet$ let $S=s_1\cup l\cup s_2\subset X$ be the closed curve composed of these three  (see the
 figure on page 175 in [G\&L 1983]).

Then let $[x_1,\delta]\subset s_1$  and $[x_2,\delta]\subset s_2$ be the  subsegments of length $
\delta$ attached to the $x_1\in s_1$  and    $x_2\in s_2$  ends of the segments $s_1$ and $s_2$  for  $\delta=dist (x_1,x_2)/3$.

Then the distances between these two subsegments, as well as between the pair  of the two complementary segments  in $S$,    that are $l$ and that the part of $S$ within distance  $\leq r-\delta$
 from $  a_0$, 
are, clearly, both $\geq \delta$;  thus, according to the above cubical example,  
$$fil.rad_\mathbb Q( \tilde S) \geq \frac{1}{2}\delta \geq  \frac {diam(L)}{6}.$$ 
QED.

 {\it Remarks.}  The above  is borrowed from  the proof of  corollary  10.11 in [G\&L 1980] claiming a {\it bound on Uryson 1-width  of 3-manifolds $X$ with $Sc(X)\geq \sigma>0$}.   
 
Since the inequality $fil.rad \leq const/\sqrt\sigma$  for 3-manifolds $X$ is crucial
for our proof of non-asphericity of 5-manifolds with positive scalar curvature and since 
the  argument following,  10.11 in [G\&L 1980] contained a flaw, we furnish   all   details needed for the proof of this inequality   here
and  in section 6.

(b) There are  similar bounds  on the diameters of connected components of 
the  "annuli" $\phi^{-1}[r,R]\subset X$, $R>r$, where, moreover, instead of the distance function to a points $x_0\in X$, one may use distance to a connected subset $X_0\subset X$.
  
  \vspace {1mm}

4.E. {\it Corollary.} {\sf If all  closed curves  $Y\subset X$ admit  finite covering 
 $\hat Y\to Y$,  such that
 $$fil.rad(\hat Y)\leq \delta,$$
  then $X$ admits a continuous map onto 
a  $1$-dimensional simplicial space, i.e. a graph $\Gamma$,
$$f: X\to \Gamma,$$
 such that the diameters  of the 
pullbacks $f^{-1}(\gamma)\subset  X$, $\gamma\in \Gamma$, are bounded by 
 $6\delta$, that is
 $$width_1(X)\leq 6\delta.$$}
 
In fact, such  a  $\Gamma$ is obtained as the  quotient of $X$ by the   relation, where two points are declared equivalent if they
 lie in the same connected  component of a level set of  some distance function on $X$, as it is explained in 
argument following corollary 10.11 in [G\&L 1983].)\vspace {1mm}

{\it Question.} {\sf Let $X$ be a path metric space, where all closed  curves $S$ {\it homologous to zero} bounds within distance  $d$ from $S$. Does then $X$ admit a map to a graph, such that the pullbacks of all points have diameters $\leq const\cdot d$?}

\vspace {1mm}

A metric space $X$ is called {\it uniformly contractible} or {\it geometrically controlled contractible}\footnote{This is called "geometrically contractible" in [M.G. 1983].} if there exists a  contractibility control function $R(r)=R_X(r)\geq r$, such that all r-balls  $B_x(r)\subset X$ are contractible in the concentric $R$-balls
$B_x(R)$, for  $R\geq R(r)$ and all $x\in X$.

A metric space  $X$ is called  {\it uniformly  acyclic} if the inclusion
homomorphisms 
$$H_i( B_x(r)) \to  H_i( B_x(R))$$
vanish for $i=1,2,...$, where the corresponding $R(r)$ is called {\it the acyclicty control function.}

Similarly, one defines  {\it uniform $\mathbb Q$-acyclicty,} where the 
homomorphisms 

$H_i( B_x(r);\mathbb Q) \to  H_i( B_x(R);\mathbb Q)$ are required to vanish.

Our main examples of such spaces come from  covering of compact spaces via the following obvious

{\it Proposition.}  {\sf If a $ X$ is a contractible, acyclic or rationally acyclic manifold; such that the action of the  isometry group is {\it cocompact} on $X$,  then  $ X$ 
is {\it uniformly} contractible or, respectively, {\it uniformly} acyclic or {\it uniformly} rationally acyclic.}

In particular, 

{\it universal coverings of compact  aspherical manifolds are uniformly contractible}.

\vspace
 {1mm}

4.F.  {\it Bounds on Relative Filling Radii of Cycles in Uniformly Acyclic Spaces.}\vspace {1mm}

{\sf If $X$ is uniformly acyclic, then the (relative) filling radii of all  submanifolds (and subpsedomanifolds)  in $X$ are bounded in terms of their absolute filling radii.}

\vspace {1mm}

More generally, let  $Y$ be a closed orientable  manifold (or psedomanifold) of dimension $n$  endowed with a metric and denote  $r=fil.rad(Y).$

Let $X$ be a uniformly acyclic space and let 
$\phi:Y\to X$ be a distance non-increasing map.\vspace {1mm}

{\sl Then the image of the fundamental cycle $[Y]\in H_n(Y)$ in $X$ is homologous to zero
in the $R^\ast(r)$-neighbourhood of the image $\phi(Y)\subset X$, where the function $R^\ast(r)$ depends only on the  acyclicty control function $R_X(r)$ and $n=dim(Y)$.}\vspace {1mm}

{\it Proof.}  To clarify, let first  $X$ be uniformly contractible with control function $R(r)$ and extend 
the map $\phi: Y\to  X$  to a map $\Phi : Y_+ 
\to X$ for   a metric extension $Y_+\supset Y$ of $Y$, where $Y$ is homologous to zero and where 
$dist(y_+, Y)\leq r$  for all $y _+,\in Y_+$.

Assume without loss of generality that $ Y_+ $ is a polyhedral space and divide it into $\varepsilon$-small simplices $\Delta$.

  Start with the map $\Phi_0$ on the 0-skeleton, of $Y_+$ obtained by sending all vertices from $Y_+$ to the  {\it nearest points} in $Y\subset Y_+$ and 
  then applying  $  \phi$.
  
  Next, extend $\Phi_0$  to   a map $\Phi_1$ from  the 1-skeleton of $Y_+$ to $X$, extend this to the 2-skeleton etc.

 At every stage, an   extension of $\Phi_{i-1}$  from the boundary of an 
   $i$-simplex  $\delta\subset Y_+$
   is possible in the $R(2r_{i-1})$-neighbourhood of the $\Phi_{i-1}$-image
   of $\partial \delta$  for $r_{i-1}$ denoting the   diameter of this image. 
 Thus, eventually,  the extension will take place in the  $R^\ast(r)$-neighbourhood of $\phi(X)$, where this $R^\ast$     is the value of the $(n+1)$-th iteration of  the function $2R(2r)$. 

Now, in the general uniformly acyclic case, instead of extending  maps from   $\partial \Delta^i$ to $\Delta^i$ we fill in $\partial \Delta^i$ at every stage by an $i$-chain  and thus eventually  fill-in  all of $ \phi_\ast[X]$ in 
the  $R^\ast$ neighbourhood of $\phi( Y) \subset X$. QED.

{\it Remark.} The above argument also yields a similar inequality for rational filling radii in uniformly $\mathbb Q$-acyclic spaces.

\vspace {2mm}

It seems, the universal  coverings $\tilde {\underline X}$  of  closed  aspherical $m$-manifolds %$\underline X$ 
in all {\it known}  examples
 contain $i$-cycles  $\underline Y_d\subset \tilde {\underline X}$, for all  $i=0,1,2,   m-1$,  with {\it arbitrary large filling radii,}
$$fil.rad (\underline Y_d\subset \tilde {\underline X})\geq  d\to \infty,$$
i.e. these $\underline Y_d$ {\it don't bound} in their $d$-neighbourhoods in $\tilde {\underline X}$,
that is equivalent to the existence of $(m-i-1)$-cycles  $\underline Y^\diamonddot \subset \tilde {\underline X}$ which have {\it non-zero} linking numbers with $\underline Y_d$  and 
$$dist(\underline Y_d, \underline Y^\diamonddot)\geq d.$$

The existence of such cycles is  obvious for $i=0, m-1$, quite easy  for   $i=1$, while the case $i=m-2$  follows by the {\it Alexander duality}  from that for $i=1$ as it  is shown in the lemma below.

On the other hand, it seems to be {\it unknown} for $n\geq 5$ 

{\sf if the universal coverings of  all compact aspherical n-manifolds contain $2$-cycles

 with arbitrarily large filling radii.}

\vspace{1mm}

4.G.{\it Codim 2  Linking Lemma}.  {\sf  Let $\tilde {\underline X}$ be a complete uniformly acyclic Riemannian manifold
of dimension $m \geq 2$ (e.g the universal covering of an aspherical manifold). 
Then 

{\it for all $d>0$ there exist a    pair of domains with smooth boundaries,  $\underline X=\underline X_d\subset \tilde {\underline X}$ and 
$\underline X_0 \subset \underline X$, such that
 
$\bullet$ $dist(\underline X_0, \partial  \underline X)\geq d$;

 $\bullet$ there exists a  homology class $\underline h_0\in H_{m-2}(\underline X_0)$ such that its image  
$ \underline h\in H_{m-2}(\underline X)$ under the inclusion homomorphism for $  \underline X_0 \hookrightarrow \underline X$ is indivisible\footnote{ That is $\underline h\neq k\underline h'$ unless $k=\pm1$.  }, hence non-zero and non-torsion.} }
  %inclusion homomorphism 
 %$$H_{m-2}(\underline X_0)\to H_{m-2}(\underline X)$$
%doesn't vanish.}

 \vspace {1mm}

{\it Proof}. Since $\tilde {\underline X}$ is uniformly  $(m-1)$-acyclic, i.e. all $(m-1)$-cycles with diameter $r$ bound in their $R$-neighbourhood for some   control function  $R=R_{\underline X,m-1}(r)$, the manifold $\tilde {\underline X}$ is also {\it uniformly connected at infinity} for $m\geq 2$.

 It follows that, for all $r>0$,

 {\sl there exists  a distance minimizing geodesic segment $\gamma\subset \tilde {\underline X}$ of length $\geq 100r$  and a curve  $C\subset  \tilde {\underline X}$ that {\it joins the ends $x_1$ and $x_2$} of $\gamma$ and  {\it doesn't intersect}  the $r$-ball around the central point $x\in \gamma$.}

Let $F:X\to \mathbb R^2_+$ be the map defined by the pair $(f_1, f_2)$  of the distance functions
$f_1(x)=dist (x, x_1)$  and  $f_2(x)=dist (x, \gamma)$.

Let   $D^2\subset  \mathbb R^2_+$ be the disc   with   center at 
the point $(r/2,r/2)\in \mathbb R^2_+$ and  and  with   radius $r/3$.

Let  
$$\Phi: X\to D^2$$ 
be the composition of $F$ with the normal   projection (retraction) $\mathbb R^2_+ \to  D^2$.

Since the  $f$-image of the closed curve
$S=\gamma\cup C$ doesn't  intersect  $D^2\subset  \mathbb R^2_+$,
the map $\Phi$ sends   $S$ to the boundary circle $\partial D^2$,  where the resulting  map $S\to \partial D^2$ has 
 {\it degree one} (or minus one depending on how we orient the two circles), because this map is {\it locally one-to-one} over the point 
$(r/2, r/6)\in \partial D^2$.

It follows that

{\sl the pullback $ Z_r= F^{-1} (r/2,r/2)\subset   \tilde  {\underline X}$} or rather,  this  pullback for a generic smooth approximation of $F$,  {\it has linking number one with $S$.}

Since $X$ is acyclic, the (transversal)  intersections of $Y_0$ with  (smooth approximations of ) the boundaries  $R$-balls around   $x$, 
$$Z_r\cap \partial B_x(R)$$
are {\it homologous to zero in the complements to the balls}  $B_x(\rho) \subset B_x(R),$ 
where 
$$\mbox {$\rho \to \infty$ for $R\to \infty$.}\footnote{Acyclicty (e.g. contractibility) of $\underline X$ implies that  the inclusions  homomorphisms
$$H_i(\partial B_x(R))\to H_i(\underline X\setminus B_x(r)),\mbox {  }i=1, 2,...,m-2,   $$   
  vanish for   
 $\rho=\rho(R)=\rho_{\underline X,x}(R)\underset {R\to \infty}\to\infty.$}$$

Let  $\underline Y_r\subset  \tilde {\underline  X}\setminus  B_x(R)$   be the $(m-2)$-cycle  obtained by attaching   
\vspace{-2mm}$$\mbox{a $(m-2)$-chain $Z' \subset \tilde {\underline X}\setminus B_x(\rho)$  with $\partial Z'= -\partial Z_r$}$$\vspace{-5mm} 

\hspace{-6mm}to $Z_r$.

%In fact, acyclicty (e.g. contractibility) of $X$ implies that  the inclusions  homomorphisms
%%$$H_i(\partial B_x(R)\to H_1(X\setminus B_x(r)$$   
%vanish for  $i=1,2,...m-2$, where  
% $$\rho=\rho(R)=\rho_{X,x}(R)\underset {R\to \infty}\to\infty.$$ 

Finally, let $ \underline X_0\subset \tilde {\underline X}$ be  a small\footnote{"Small" means that it is contained in the $\varepsilon$ neighbourhood, where, eventually, $\varepsilon\to 0$.} 
  smooth  neighbourhood  of  $\underline Y_r$, where $R_+\geq R$, and let $\underline X =\underline X_d\supset \tilde{ \underline X}$ be a  smoothed  $d$-neighbourhood of $\underline X_0$.

Since the cycle $\underline Y_r$ linked  with the curve  $S$ in in $\tilde {\underline X}$   with linking number one, 
the homology class $\underline h\in H_{m-2} (\underline X)$  represented by
 $\underline Y_r  \subset \underline X =  \underline X_d$   is indivisible, where   $d\to \infty$ for  $r\to \infty$ and provided that  $R$ as well as $\rho$ are much larger than $r$. QED.

\section {Non-compact Version of the  Chodosh-Li Theorem}

 {\sf Let  $\tilde X$ be  a complete orientable Riemannin $4$-manifold, which admits a proper  {\it distance non-increasing} map to a   complete {\it uniformly acyclic}  Riemannian 4-manifold  $\tilde {\underline X}$
 $$f: \tilde X \to \tilde {\underline X}$$
 and let the degree of this map be non-zero.
 Then the lower bound $\sigma$ on the scalar curvature of $\tilde X$ is {\it non-positive},
 $$\sigma=\inf_{x\in \tilde X} Sc(\tilde X,x)\leq 0.$$
 
 Consequently,   {\sl compact aspherical 4-manifolds admit no metrics with $Sc>0$}.}

 \vspace {1mm}

{\it Proof. } Let $(\underline X_d, \underline X_0)$,  
$$\tilde {\underline X}\supset \underline X_d\supset \underline X_0,$$ 
 be a pair delivered by the above codim 2 linking  Lemma and let 
 $(X, X_0)$  be the $f$-pullback of this pair,
 $$X= X_d=f^{-1} (\underline X_d)\subset \tilde X\mbox { and } X_0=f^{-1} (\underline X_0)\subset \tilde X.$$

Since  $k=deg (f)\neq 0$, 
there exists a homology  class 
$h_0\in H_{m-2}(X_0)$, namely  (Gysin's) pullback $Y_d=f^{-1}( \underline Y_d)$ of   the cycle $\underline Y_d$), such that the class 
$f_\ast (h_0)=  \underline h_0\in H_{m-2}(\underline X_0)$ doesn't vanish in $\underline X\supset \underline X_0$. 

Since $f$ is distance non-increasing, 
$$dist(X_0, \partial X)\geq dist (\underline X_0, \partial \underline X_d)$$
and, also,

 {\sf $f$ doesn't increase the filling radii of  cycles
in $\tilde X$.}

It follows that 

\hspace {0mm}{\it the filling radii of all cycles in $X_0$ in the class $h_0$, i.e.  homologous to  $Y_d$,  

are greater  than or equal to  $d$.}

On the other hand,  by Richard's lemma, there exists 
   a  closed surface $Y\subset X$, homologous to $Y_d$,    such that 
a warped  $\mathbb T^2$-extension $Y^\ast$ of  $Y$ has 
$Sc(Y^\ast)\geq \sigma/2$.

According to Example 1, all connected components  $Y_i$ of this $Y$ have diameters $\delta \leq \frac{4\pi}{\sqrt \sigma}$
and  since the map is distance non-increasing,  the $f$-images  $\underline  Y_i\subset \underline X_0$  of these components  have diameters $\underline \delta\leq \delta \leq \frac{4\pi}{\sqrt \sigma}$ as well.

Since the ambient manifold  $\tilde {\underline X}\supset \underline X_0$ is uniformly acyclic,
the filling radii of all $ \underline  Y_i$, hence of $ \underline  Y=f(Y)$, are bounded by 
 $R(\underline \delta)$, where this quantity  remains bounded for $\sigma>0$,  while $d\to \infty$;
 hence $\sigma=0$. QED.

\vspace {1mm}

{\it Remark  on Non-complete Manifolds.}   All of the above remains true for compact manifolds with boundaries, in-so-far we operate sufficiently far from these boundaries.

It follows, for instance, that  if $ \underline X$ is complete uniformly acyclic $4$-manifold, then for all $\sigma>0$, there exists $R=R_{\underline X}(\sigma)$, such that \vspace {1mm}

{\sl no (possibly non-complete)
$4$-manifold   $X$  with $Sc (X)\geq \sigma$  admit a proper distance decreasing map of positive degree 
onto an (open)  $R$-ball in $ \underline X$.}\vspace {1mm}

{\it Quadratic Decay Corollary.} {\sf Let $X$ and $ \underline X$ be complete orientable  Riemannian 4-manifolds, where
$Sc(X)>0$ and $ \underline X$, is uniformly acyclic.

Let $f:X\to \underline X$ be a proper smooth  map of non-zero degree.}

Then 

{\it the ratio $\frac {Sc(X,x)}{||df(x)||^2},$
where $df:T(X)\to T(\underline  X)$  stands for  the differential of $f$ decays at least quadratically on $X$, that is 
the quantity  
$$\cdot \inf_{x\in B_{x_0}(R)} \frac {R^2\cdot Sc(X,x)}{||df(x)||^2}$$ 
remains bounded for $R\to \infty$, where  $ B_{x_0}(R)\subset X$ denotes the R-ball around a fixed
 point $x_0 \in X$.}

\section {Bounds on Uryson's Widths   and Filling Radii   of 3-Manifolds  with  Stabilized Scalar Curvatures $\geq \sigma>0$}

\vspace {1mm}

In order to extend the above  4D-argument to 5-dimensional manifolds we need to show that 

{\sl closed  3-manifolds $X$,
which admit warped $\mathbb T^2$-extensions $X^\ast$ with $Sc(X^\ast)\geq \sigma>0$, have their (absolute)  filling radii  bounded by 
$$fil.rad (X) \leq const_3 \sqrt \sigma.$$} 
This is proven in this section below, where, in fact, we  show that that  {\it Uryson's 2 -width} of such an  $X$ is bounded  by $const'_3 \sqrt \sigma.$
\footnote{
  It is sufficient to have such a bound for  a {\it finite, or even  amenable,  covering}  $\hat X$ of $X$, which may be useful in some cases.}

Namely we prove the following 

\vspace {1mm}

6.A. {\it 3d-Covering Theorem.}   {\sf Let a  closed orientable  Riemannin $3$-manifold   $X$ admit a warped 
 $\mathbb T^N$-extension 
 $X^\ast$, such that  
$$Sc(X^\ast)\geq \sigma.$$
{\it Then  there exists  a finite 
collection of  closed subsets $ V_i$   that cover $  X$,  
$$X=\bigcup_i V_i,$$
such that 

$\bullet $ the multiplicity of this covering is $\leq 3$, i.e. no four of these  subsets intersect; 

 $\bullet $  the diameters of all these subsets are bounded by a constant $D$ depending only on $\sigma$, namely,
 $$diam(V_i)\leq D\leq \frac {12\pi}{\sqrt \sigma}+\varepsilon.$$}}
 
6.B.  {\sc Filling Corollary.}
 $$fil.rad(X)\leq  \frac {6\pi}{\sqrt \sigma}+\varepsilon.$$

%\hspace {30mm} {\sc Remarks and Corollaries.} \vspace {1mm}
 
 %{\it Remark.}   It is implied by  corollary  10.11 in [G\&L 1980]  that  if $Sc(X)\geq \sigma$, then $X$  admits a covering by closed subsets $V_i$ with multiplicity two (rather than three), where  
 %$diam(V_i)\leq \frac {12\pi}{\sqrt \sigma}.$ But the argument suggested in [G\&L 1983] delivers such a covering only for  manifolds $X$ with $H_1(X;\mathbb Q)=0$.
 
%COROLLARY 10.11

 %Mikhael Gromov and H. Blaine Lawson, Jr., Spin and scalar curvature in the presence
%of a fundamental group. I, Ann. of Math. (2) 111 (1980), no. 2, 209{230. MR 569070
%[6] , Positive scalar curvature and the Dirac operator on complete Riemann-
%ian manifolds, Inst. Hautes Etudes Sci. Publ. Math. (1983), no. 58, 83{196 (1984).\S12

\vspace {1mm}

{\it Proof.}  We shall  obtain the required covering of $X$ by firstly  cutting $X$  into    submanifolds $X_j$  with  {\it boundaries},  such that  the homology  groups $H_1(X_j)$  are pure torsion 
  and  then applying   the   argument following  corollary 10.11 in   [G\&L 1983] to these $X_j$.

 To clarify the idea, let us start with the case, where the manifold $X$ itself has $Sc(X)\geq \sigma$  and recall that  
the fundamental group of such an $X$ is a free product of finite and infinite cyclic  groups, see e.g. section 8 in [G\&L 1983].

 Let $\{\Sigma^2_k\}$  be a {\it maximal}  collection of  {\it disjoint } connected  stable minimal surfaces in $X$, where  "maximal" means that every  such connected  surface in the complement $ X\setminus \cup_k \Sigma^2_k$ is isotopic to one of  $\Sigma^2_k$.
 
 Since all these $ \Sigma^2_k$ are topological spheres by the Schoen-Yau argument from [S\&Y 1979] and since  all 2-dimensional homology classes in manifolds with mean convex, e.g. minimal, boundaries 
 are realizable  by stable minimal  surfaces,  all connected components $X_j$   of $X\setminus \cup_k \Sigma^2_k$  have their 1-dimensional homology {\it pure torsion}, i.e. $H_1(X_i;\mathbb Q)=0$.
 In fact, the fundamental groups of $X_i$,  are free products of finite groups.\footnote{If one feels uncomfortable with torsion, one may take  a finite covering  $\hat X\to X$ with a {\it free fundamental group}  $\pi_1(\hat X)$ and apply what follows to $\hat X$. This  will  result in  a lower bound  on the {\it rational} filling radius of $X$, i.e. defined with  the  {\it rational homology}, $H_3(X;\mathbb Q)$ rather than $H_3(X)=H_3(X;\mathbb Z)$. 
 
 This causes no problem,  
 since the class $\underline h$ in the  codim 2 linking lemma in section 4 is indivisible; hence, remains non-zero after tensoring with $\mathbb Q$.}

 Also, since the boundaries of all connected components $X_j$  of  $  X\setminus \cup_k\Sigma^2_k$  are mean convex,  a multiple of  every closed  connected curve $S$ in $X_j$
  bounds a stable minimal surface $Y\subset  X_i$, which, according to theorem 10.7 in [G\&L 1983], lies within distance $\leq 2\pi/\sqrt \sigma$ from $S$.\footnote {The argument indicated for the proof of example 1 in section 1 yields the bound $dist(y, \partial S)\leq \sqrt 2/\sqrt\sigma$, $y\in Y$.}
 
 \vspace {1mm}
 
Now it follows from 4.E  that each of the above 3-manifolds  $ X_j$ with $\Sigma$-boundaries  admits a map 
 to a finite  $1$-dimensional simplicial space, i.e. a graph $\Gamma_j$,
$$f_j: X_j\to \Gamma_j,$$
 such that the diameters  of the 
pullbacks $f^{-1}(\gamma)\subset  X_j$, $\gamma\in \Gamma_j$ are bounded by 
 $\frac {12\pi}{\sqrt \sigma}.$

Cover  the graphs $\Gamma_j$   by $\epsilon$-small subsets  
 $\Gamma_{jl}$  with multiplicities 2 and then  cover $  X$  
 
$\bullet$ by  {\sf $\epsilon$-neighbourhoods $U_{k\epsilon}\subset X$ of $\Sigma_k^2$ 

and 

$\bullet$ by  the closures of the differences 
$$\bar X_{ j, k,l,\epsilon}= f^{-1}(\Gamma_{j,l})\setminus   U_{k\epsilon}  \subset  
 X_j \subset X \mbox { for all   $j,k$ and $l$}.$$}
  The resulting  covering  of $X$, call it $\{V_{i, \epsilon}\}$, clearly, has multiplicity $\leq 3$, where  the diameters of  $U_{k\epsilon}$ are bounded by 
 $2\pi/\sqrt\sigma+2\epsilon$ (Example 1) and the pullbacks $f^{-1}(\Gamma_{j,l})$ satisfy 
  $$diam (f^{-1}(\Gamma_{j,l}) \leq
 12\pi\sqrt\sigma+ \varepsilon,\mbox { where }\varepsilon =\varepsilon(\epsilon)\underset{\epsilon \to 0}\to 0.$$
 This furnishes the proof in the case $Sc(X)\geq \sigma$.
 
\vspace {1mm}

\vspace {1mm}

Now, let us turn to the general case, where

\hspace {10mm} {\sl   a warped $\mathbb T^N$-torical extension $X^\ast$  of $X$ has  $Sc(X^\ast)\geq \sigma>0.$}

 It is convenient at this stage to work with the quotient  manifold  $\underline X^\ast=X^\ast /\mathbb Z^N$ (see section 1), where, observe,  
 stable minimal $\mathbb T^N$-invariant  hypersurfaces in  $X^\ast/\mathbb Z^N$ correspond to
 local minima of the following functional  on closed surfaces $Y\subset X$,
 $$Y\mapsto \int_Y vol(\mathbb T^N_y)dy,$$
 where $\mathbb T^N_y$ are  the torical fibers $\mathbb R_y^N/\mathbb Z^n$ (see section 1).
 
 Thus the above existence theorems for stable minimal surfaces in   $X=(X, g)$  apply to $X$ where  the  original metric
 $g$  is  modified by the conformal factor, 
 $$\mbox {$g\leadsto \psi(x) g(x)$ for $\psi(x)=vol(\mathbb T^N_x)$, $x\in X$}.$$

 On the other hand, the necessary bounds  on distances in these surfaces with respect to the original metric $g$  follow by the argument used in example 1. 
 
 This accomplishes the proof of the  3d covering theorem.
 \vspace {1mm}

{\it 3d Geometric Decomposition Theorem*.} {\sf Let $X$ be  a compact orientable Riemannian 3-manifold with $Sc(x)> \sigma>0$. Then there exists a a collection of disjoint smooth embedded spheres $\Sigma_\mu\subset X$ with the following properties.}

 (i) {\sl The diameters of all $\Sigma_\mu$ in the induced Riemannin metrics  satisfy
 $$diam(\Sigma_\mu)\leq \frac {3\pi\sqrt 2}{\sqrt\sigma}, $$}

 (ii) {\sl the diameters of the connected components $X_\nu$ of the complement $X\setminus\cup_\mu \Sigma_\mu$ satisfy
  $$diam(X\nu)\leq \frac { (2\sqrt 6+3\sqrt 2)\pi}{\sqrt\sigma}, $$}

 (iii) {\sl all  $X_\nu$ have finite fundamental groups; hence, they are  
 diffeomorphic to the spherical 3-forms  with finitely many punctures.}
 
 \vspace{1mm}
 
 {\it Outline of the Proof.} We start with  cutting $X$ by minimal spheres 
 to pieces $X_j$  as earlier which  satisfy (iii) but not necessarily (ii). 
 
 Next we cut   each $X_j$
 by discs  with punctures, call them $ \Upsilon _{jk}\subset X_j$, into pieces with  small  diameters by using a version of theorem 2 for manifolds with boundary.  
 
 Then we move the boundaries  $\partial\Upsilon _{jk}\subset \partial X_j$  inside $X_j$ without creating new intersections.
 
 (Such a "movement" only  insignificantly increases the diameters of $\Upsilon _{jk}$, but it may uncontrollably enlarge their areas and it  it remains  unclear if one can have the  areas of the spheres $\Sigma_\mu$    bounded by  $const/\sigma.$)

 Since we don't use this in the sequel we leave the actual proof to the reader.\vspace {1mm}
 
 {\it Problem 1}. {\sf Can one  use the Hamilton Ricci flow to  {\it canonically} decompose 
  complete 3-manifolds $X$   with $Sc(X)\geq \sigma \geq 0$,  and more generally,
  for  $Sc(X^\ast)\geq \sigma>0$, into $X_\nu$-like pieces?
 
 (This would imply, for instance that compact aspherical Riemannin manifolds of all dimensions $n\geq 3$ admit  no 3-dimensional foliations with leaves with positive scalar curvatures.)
 
  {\it Problem 2}. Is there a counterpart to the  3d-covering and decomposition theorems for 4-manifolds?}

 %[SY(incompressible) 1979] R.Schoen, S.T Yau {\sl  Existence of incompressible minimal surfaces and the topology of three dimensional
%$manifolds of non-negative scalar curvature,} Ann. of Math. 110 (1979), 127-142.\vspace{1mm}  \vspace{1mm}

\section {Proof of the Main Theorem and an Alternative Proof of  Chodosh-Li theorem}
 
Let us reformulate the theorem stated in the summary  in the  the $\mathbb T^\ast$-stable quadratic decay form  as in  the corollary to the non-compact 
Chodosh-Li theorem.\vspace{1mm}

{\it Main Theorem.}  {\sf Let $\tilde X$ and $ \tilde {\underline X}$ be complete 
orientable  Riemannian 5-manifolds,  let $\tilde X^\ast$ be a warped $\mathbb T^N$-extension of $\tilde X$,
such that $Sc(\tilde X^\ast,x)>0$, $x\in \tilde X$,\footnote{Recall that the function $Sc(\tilde X^\ast)$, on $\tilde X^\ast$  being $\mathbb R^N$-invariant, is regarded as a function on $\tilde X$. }
and let $ \tilde {\underline X}$ be uniformly acyclic.

Let $f:X\to \underline X$ be a proper smooth  map of non-zero degree.}

Then 

{\it the ratio $\frac {Sc(\tilde X^\ast,x)}{||df(x)||^2},$
where $df:T(\tilde X)\to T(\tilde  {\underline  X})$  stands for  the differential of $f$, decays at least quadratically on $\tilde X$, that is 
the quantity  
$$\cdot \inf_{x\in B_{x_0}(R)} \frac {R^2\cdot Sc(\tilde X^\ast,x)}{||df(x)||^2}$$ 
remains bounded for $R\to \infty$, where  $ B_{x_0}(R)\subset \tilde X$ denotes the R-ball around a fixed
 point $x_0 \in \tilde X$.}\vspace {1mm}

{\sc Main Corollary}. {\sf let  $\tilde X$ be the universal covering of a compact manifold $X$, let 
$N=0$, i.e.  $\tilde X^\ast$, let    also  $\tilde {\underline X}=\tilde X$ and let the map $f$ be the identity. Then the scalar curvature of $\tilde X$ can't be everywhere positive.  This means  that  compact  aspherical 5-manifolds   can't have 
$Sc > 0$.}
 
(A posteriori, if $Sc(\tilde X)\geq 0$, then $\tilde X$ is isometric to the Euclidean space $\mathbb R^5$.)\vspace {1mm}

{\it Proof.}  We argue  as in the proof of the non-compact Chodosh-Li theorem, in section 5 except that instead of a surface $Y \subset \tilde X$   with large filling radius, obtained via  a curve  $S\subset \tilde{\underline X}$, we have (again by  codim 2  linking Lemma 4.G) a closed orientable  $3$-manifold, now called $X\subset \tilde X$, with $Sc (X^\ast) \geq \frac {3}{5}\sigma$, the $f$-image of (the fundamental class of) which in $\tilde{\underline X}$ is  linked to such an $S$. 
 
 Besides, instead of directly filling-in $f(X)\subset \tilde{\underline X}$  as we were doing it with $f(Y)$, which had diameter bounded by $const/\sqrt\sigma$,  we apply 
 the bound on the absolute filling  radius of $X$ in terms of $\sigma =\inf Sc(X^\ast)$ 
 (filling corollary 6.B from  the previous section) combined with the bound  4.F  on the relative filling radii of cycles  in uniformly acyclic spaces.

      Thus  we  conclude to 	the vanishing  of the fundamental  class $[X]$ of $X$ in its 
 $d$-neighbourhood  in $\tilde X$. QED.

\vspace {1mm}

{\it Codim 1 Proof of the Chodosh-Li Theorem.} Let $\tilde X$ be a complete non-compact Riemannian manifold of  with $Sc(\tilde X)\geq \sigma>0$.

If $m=dim(\tilde X)\leq 8$ then

{\it the manifold $\tilde X$ can be exhausted  by compact domains with smooth boundaries, 
$$X_1\subset X_2\subset... \subset X_i \subset .... , \bigcup_iX_i =X,$$
such that the boundaries  $Y_i=\partial X_i$ of all $X_i$ admit warped $\mathbb T^1$-extensions $Y_i^\ast$  with $Sc (Y_i^\ast)\geq \frac{m-2}{m} \sigma$.} 
 
 In fact, by theorem 2, one  may even  choose  $X_i$, such that 
 $$dist (\partial X_i, X_{i-1})\leq\frac { 2\pi m(m-1)}{\sqrt \sigma}.$$
 
 On the other hand, if $\tilde X$ is uniformly acyclic, 
 then, obviously,  the filling radii of the fundamental classes of $Y_i$ tend to infinity,
  which, for $n=m-1=3$ implies that $\sigma\leq 0$ by the filling corollary 6.B to 3d covering  theorem in the previous section. QED. 
 
\vspace {1mm}

 \section {Filling radii of Submanifolds in manifolds $X$ with $Sc(X) \geq \sigma>0$}
 
 The geometric inequalities used in the above arguments   can be regarded as generalizations of the bounds on the filling radii of closed curves in 3-manifolds.
 
 In fact, our proofs of these inequalities applied to minimizing submanifolds with given boundaries $Y\subset X$  imply the following.
 
{\it Proposition. }  {\sf  Let $X$ be a compact orientable  Riemannin  manifold with $Sc(X)\geq \sigma>0$ of dimension $m\leq 5$ and let $Y\subset X$  be codimension 2  submanifold that is homologous to zero.
Let $Y$ be endowed with  the metric, where the distance equals to that measured in $X$.}
$$dist_Y(y_1,y_2)=dist_X(y_1,y_2),\mbox {   }  y_i,y_2\in Y\subset X.\footnote {This $dist_Y$ may be significantly smaller than the distance associated  to the Riemannian metric in $Y$ induced from $X$.}$$
%Then {\it the  filling radii of  orientable codimension 2  submanifolds $Y\subset X$  with the distance functions coming from $X$, 

Then {\it the filling radius of $Y$ is  
   bounded by $const/\sqrt \sigma$  for $const \leq 100$. 
More generally, if $X$ is a manifold with boundary, where  $dist(Y, \partial X)\geq 100/\sqrt\sigma$, then the filling radius of $Y$ satisfies the same  bound, provided  the homology class $[Y]\in H_{m-2}(X) $ is contained in the image of the inclusion homomorphism $H_{m-2)}(\partial X) \to H_{m-2}(X).$

Furthermore if $m=dim(X)\leq 4$,  similar inequalities hold for codimension one submanifolds $Y\subset X$.}

 However, unlike the 3-dimensional case, these inequalities don't yield bounds on $fil.rad(X)$, where the difficulty resides with "elementary metric lemma" 4.D, which has no apparent higher dimensional counterpart.

 \section {Perspectives for $m>5$}
 The main theorem, combined with theorems 1 and 2, implies non-existence of metrics with $Sc\geq \sigma >0$ on special manifolds of dimensions $m>5$, at least for $m\leq 8$.
  
   For instance, 
 \vspace{1mm}
 
 {\it if a complete Riemannian aspherical manifold $X$ is homeomorphic 
 to $X_0\times \mathbb R^2\times \mathbb T^1$, where $X_0$ is a compact 
 5-manifold, then $ \inf_xSc(X,x)\leq 0$.}  \vspace{1mm}
 
 Also, as it   is observed in [J.W. 2019],
 
 {\it  images of the rational  fundamental classes of compact $3$-manifolds  with $Sc>0$ 
 vanish under continuous maps into  aspherical spaces.\footnote{This  is because these manifolds are   connected sums of $S^2\times S^1$ and spherical space forms.  }}
 (This is unknown for $m$-manifolds for $m\geq 4$.\footnote{If a compact $m$-manifold $X$  admits,   say  piecewise linear, map  $\phi$ to a $(m-1)$-dimensional polyhedral space $P$,  such that  the images of the fundamental groups of the connected components of the  pullbacks $\phi^{-1}(p)\subset X$  in the fundamental group of $X$ under the  inclusion homomorphisms,
 $$\pi_1(\phi^{-1}(p))\to\pi_1(X),$$
 are {\it finite} for all $p\in P$, then
 all continuous maps from $X$ to aspherical spaces send the fundamental rational homology class $[X]_\mathbb Q$ to zero. But it is  unlikely(?)  that  all compact Riemannian manifolds with $Sc>0$  
 admit   $\phi$  with  this property  for $m=dim(X) > 3$.})

 Thus for instance,    
 
 {\it if a closed aspherical manifold $X$ of dimension $m=3,4,....,8$, admits a continuous map $f$  to  the torus $\mathbb T^{m-3}$,  such that the  (homology class of the ) point-pullback   $f^{-1}(t)\subset X$,  $t\in \mathbb T^{m-3}$, doesn't vanish in $H_{m-3}(X;\mathbb Q)$,  then $X$  carries no metric with $Sc>0$}. (This is unknown for such maps  $f:X\to \mathbb T^{m-n}$ for $m\geq 6$ and $ n\geq 4$.)

 \vspace{1mm} 
 
 In general, there are two problems of quite different kinds  for proving  non-existence of  metrics with $Sc>0$ on aspherical  manifolds of high dimensions  that we formulate below  for $m=6$ in the form of two {\it \large conjectures}.

  A. {\it "The Universal Coverings of Compact acyclic Manifolds are Uniformly Acyclic  in Codimension One":} {\sf  Let $\tilde X$ be the universal covering of a  compact aspherical $6$-manifold. Then, for every $d>0$ there exits a (Riemannian) band $Z\subset  \tilde X$
  bounded by a pair of (infinite) hypersurfaces $\partial_-$ and $\partial_+$, such that
  $width(Z)=dist (\partial_-, \partial_+)\geq d$ and such that all (infinite) hypersurfaces 
  $Y\subset Z$ which separate $\partial_-$  from $\partial_+$  admit distance   non-increasing\footnote{Here as 
  well  as everywhere earlier,   "distance        non-increasing" conditions on maps $f$,  may be replaced by   $dist (f(x_1),f(x_2)) \leq \theta(dist(x_1,x_2))$, where  $\theta(d)=\theta_f(d)$ is an arbitrary continuous (hence locally bounded)  function in $d\in [0, \infty )$.}  
    maps  of degrees $\neq 0$ to uniformly rationally acyclic 5-manifolds.\footnote{Examples  of such bands are those between pairs of  parallel hyperplanes in the Euclidean space  $\mathbb R^6$ and, more generally, between concentric horospheres in (complete simply connected)  manifolds with non-positive sectional curvatures.}

  \vspace {1mm}
  
 {\sl Validity of either  A or B would imply  that compact aspherical 6-manifolds carry no metrics with $Sc>0$.}

 B. {\it "The Fundamental Classes of Manifolds with Large Scalar Curvatures are Small"}:{ 
Compact aspherical 4-manifolds $Y$  with $Sc(Y)\geq \sigma>0$ have their filling radii bounded by
 $$fil.rad_\mathbb Q(Y)\leq \frac{const_4}{\sqrt \sigma}.$$}}

    \vspace {1mm}

   {\it Remark.} The solution of B doesn't seem sufficient for ruling out maps from $Y$ to aspherical spaces 
  with non-zero images of their fundamental classes $[Y]\in H_4(Y, \mathbb Q)$,
  since the following  following seems unsettled.

{\it Question.}   Let let $Z$ be an aspherical space, $Y$ be an orientable  $n$-dimensional manifold or pseudomanifold,
such that all covering spaces $\hat Y\to Y$  have their filling radii (defined below) bounded by
$$fill.rad(\hat Y)\leq const=const(Y),$$
  and let $f:X\to Z$ be a continuous map.
  
  Does then the rational  homology mage $f_\ast[Y]_\mathbb  Q\in H_n(Z;\mathbb Q)$ vanish?
  
  \vspace{1mm}

  {\it Definition of  filling radii of manifolds with boundaries and   non-compact manifolds.}(Compare with 4.4.C in [M.G. 1983]).  Let $Y$ be an orientable manifold or a pseudomanifold  of dimension $n$  with boundary.
  
   By definition, the inequality 
  $fil.rad(Y) \leq R$ for a given metric on $Y$ signifies that    there exists a  metric extension $Y_+\supset Y$, such that  $dist (y_+, Y)\leq R$, $y_+\in Y_+$, and such that  the
   relative fundamental class  $[Y]\in H_n(Y, \partial Y)$   vanishes under the inclusion homomorphism
    $$H_n(Y, \partial Y)\to H_n(Y_+, U_R(\partial Y)),$$
  where  $U_R(\partial Y)\subset Y_+$ denotes the $R$-neighbourhood of $\partial Y\subset Y\subset Y_+$.
  
  Then the inequality  $fil.rad(\hat Y)\leq R$ for a  non-compact (pseudo) manifold is understood as the 
  existence of an exhaustion of $\hat Y$  by compact 
  $Y_1 \subset  Y_2 .... \subset Y_i \subset ...\subset \hat Y$
  with boundaries, such that $fil.rad(Y_i)\leq R$  fo all $i$.
  
For instance,
  $fil.rad(\hat Y)\leq\frac{1}{2}   width_{n-2}(\hat Y).$
  
  (This  applies to complete $3$-manifolds $\hat Y$ with $Sc(\hat Y)\geq \sigma>0$ and 
  shows that their filling radii are bounded by $const/\sqrt\sigma.)$\vspace{1mm}
  
 {\it Exercise.} Show that that the above question has positive answer if the fundamental group of $Z$ is {\it residually finite}, or, more generally,   {\it residually amenable}.

  %fails deliver upper bounds on the {\it "absolute" filling radii} of   submanifolds $Y\subset X$,   where
 %a warped $\mathbb T^\ast$ -extension $X^\ast$ of $X$  has $Sc(X\geq \sigma>0$ in the cases
 %$dim(X)=4$ and $codim(Y)\leq 2$ and also for $dim(X)=5$ and $codim(Y)=2$, where these estimates 

 %{\it Remark.}  A relative version of this  argument, applied to  hypersurfaces $Y$ with boundaries  in  a $4$-manifold  $X$ with $Sc(X)\geq \sigma >0$ yields the following.   

 \section {References}

\hspace {3.4mm} [C\&L  2020]  O. Chodosh, C. Li {\sl Generalized soap bubbles and the topology of manifolds with positive scalar curvature},  arXiv:2008.11888v1
 \vspace {1mm} \vspace {1mm}

 [F-C\&S 1980] D. Fisher-Colbrie, R. Schoen  The structure of complete stable minimal
surfaces in 3-manifolds of non-negative scalar curvature, Comm. Pure
Appl. Math., 33 (1980) 199-211. \vspace {1mm} \vspace {1mm}

[M.G. 1983]   M.Gromov, {\sl  Filling Riemannian manifolds}, J. Diff. Geom. 18 (1983), 1-147.\vspace {1mm} \vspace {1mm}

 [M.G. 2018] M. Gromov, {\sl  Metric Inequalities with Scalar Curvature,}  arXiv:1710.04655
 \vspace {1mm} \vspace {1mm}

[M.G.2019]  M. Gromov,{ \sl Four lectures on scalar curvature}, arXiv:1908.10612
\vspace {1mm} \vspace {1mm}

[G\&L 1980]  M. Gromov, B Lawson, 
 Spin and scalar curvature in the presence
of a fundamental group. I, Ann. of Math. (2) 111 (1980), no. 2, 209-230. \vspace {1mm} \vspace {1mm}

[G\&L 1983]  M. Gromov, B Lawson, 
{\sl Positive scalar curvature and the Dirac operator on complete Riemann-
ian manifolds},  Publ. Math. de IHES , no. 58(1), 83-196. \vspace {1mm} \vspace {1mm}
 
[L.G. 2010] L. Guth, {\sl  Metaphors in systolic geometry}. Proceedings of the International
Congress of Mathematicians 2010 (ICM 2010)
\vspace {1mm} \vspace {1mm}

[L.G. 2017] L. Guth, {\sl Volumes of Balls in Riemannian Manifolds and Uryson Width},
Journal of Topology and Analysis 9.02 (2017): 195-219.\vspace {1mm} \vspace {1mm}

[M.K. 2007] M. Katz, {\sl Systolic geometry and topology}. With an appendix by J. Solomon. Mathematical Surveys and Monographs, volume 137. American Mathematical Society, 2007.\vspace {1mm}\vspace {1mm}

 [J.L.2018 ]  J. Lohkamp, {\sl Minimal Smoothings of Area Minimizing Cones},
arxiv:1810.03157

\vspace {1mm} \vspace {1mm}

[S\&Y 1979] R. Schoen, S. Yau, {\sl  Existence of incompressible minimal surfaces and the topology of three dimensional
manifolds of non-negative scalar curvature,} Ann. of Math. 110 (1979), 127-142.
\vspace {1mm} \vspace {1mm}

 [T.R. 2020]  T. Richard   {\it On the 2-systole of stretched enough positive scalar curvature metrics on $S^2 \times  S^2$}, arXiv:2007.02705

 \vspace {1mm} \vspace {1mm}

[J. W.  2019] J. Wang, {\sl Contractible 3-manifolds and positive scalar curvature}, Ph.D. thesis, 
Universit\'e Grenoble, Alpes,
2019.

 \vspace {1mm} \vspace {1mm}

[S.W. 2007] S. Wenger {\sl A short proof of Gromov's filling inequality,}
arXiv:math/0703889.

 \vspace {1mm} \vspace {1mm}

 [J.Z. 2019]  J. Zhu, {\sl  Rigidity of Area-Minimizing 2-Spheres in n-Manifolds with Positive Scalar Curvature}   	arXiv:1903.05785

 %Katz, M. Systolic geometry and topology. With an appendix by J. Solomon. Mathematical Surveys and Monographs, volume 137. American Mathematical Society, 2007.

%Gromov, M.: Filling Riemannian manifolds, J. Diff. Geom. 18 (1983), 1?147.

%Guth, Larry. ?Volumes of Balls in Riemannian Manifolds and Uryson Width.? Journal of Topology and Analysis (February 22, 2016)

%Guth, Larry. Metaphors in systolic geometry. Proceedings of the International
%Congress of Mathematicians 2010 (ICM 2010)

%Sep 9, 2019 - Uryson width and volume. Authors:Panos Papasoglu
%Geometric and Functional Analysis volume 30, pages574?587(2020)

%A short proof of Gromov's filling inequality
%Stefan Wenger 2007
%arXiv:math/0703889 
  	
\end{document}